\documentclass[10pt]{extarticle}
\usepackage{amsthm}
\usepackage{mathtools}
\usepackage{color}
\usepackage{amsfonts}

\title{A double-inductive proof of Moessner's theorem}
\author{Archy Will He}
\date{5th Feb 2016}

\begin{document}

\maketitle

\abstract{We present a proof of Moessner's theorem by double induction, using only basic rules of arithmetic. No prerequisite knowledge is assumed.}

\section{Introduction}

Suppose we have a stream of positive integers:
$$ 1, 2, 3, 4, 5 ....$$
For any natural number $n$, remove every $x \coloneqq (n+2)$-th element and form a new stream by doing partial sum on the resulted stream. Repeat the removal-then-new-stream-formation process for $n$ times with $x$ decreases by 1 each time and we would obtain a new stream
$$ 1, M(0,n), M(1,n),M(2,n),M(3,n)...$$
where
$$\forall n \in \mathbb{N}\; \forall k \in \mathbb{N}\quad M(k,n) = (k+2)^{n+2} $$
This is one way of expressing Moessner's theorem.

To get the intuition behind it, let us start off by dropping every second number in a stream of positive integers. We obtain the stream
$$1, 3, 5, 7, 9...$$
which is a stream of odd numbers. If we construct a new stream by partial summation on this stream, the resulting stream is that of square numbers
$$1, 4, 9, 16, 25...$$
Intriguing. So what if we are to begin with dropping every third number, constructing a new stream by partial summation, then dropping every second number in that stream, and constructing a new stream again by partial summation? The resulting stream would then be the stream of cube numbers
$$1, 8, 27, 64, 125...$$
Moessner conjectured \cite{moessner} that the procedure above can be generalized into obtaining the stream
$$1^n, 2^n, 3^n, 4^n, 5^n...$$
if the procedure is began by dropping every $n$-th number, and this was what would later be known as Moessner's theorem, after first proven by Perron \cite{perron}. Moessner's theorem has been sequently proven by many others such as Paasche \cite{p}, Long \cite{l}, Hinze \cite{h}, Niqui and Rutten \cite{n_r}, and Kozen and Silva \cite{k_s}.

In this article we present a new proof for Moessner's theorem by double induction, using only basic rules of arithmetic. No prerequisite knowledge is assumed. Familiarity with summation is advised for page 7.

Throughout the article the predicate $\phi_M(k,n)$ would be used to express that the theorem holds for some $k$ and $n$ i.e.,
$$\phi_M(k,n) \Leftrightarrow M(k,n) = (k+2)^{n+2}$$
and the procedure described in the first paragraph to obtain
$$1, M(0,n), M(1,n), M(2,n), M(3,n)...$$
would be refer	red to as \textit{Moessner's Sieve}.

\section{An inductive proof for \large{$\forall n \in \mathbb{N}\; \small{M(0,n)=2^{n+2}}$}}
It is easy to see that the base case (n=0) holds i.e, $M(0,0)=2^2=4 \; \land \; 1 + 3 = 4$. We now demonstrate that $\forall n \in \mathbb{N} \; \phi_M(0,n) \Rightarrow \phi_M(0,n+1)$ by showing that
$$\forall n \in \mathbb{N}\quad M(0,n+1)= M(0,n) + M(0,n)$$
For any $ n \in \mathbb{N}$, by \textit{Moessner's Sieve}, $M(0,n)$ is equivalent to the sum of $A_n\coloneqq(n+2)+1$ and $B_n$ where $B_n$ is the sum of the numbers below:
\begin{align}
\begin{split}\label{triangle1}
&f_0(0),f_0(1),f_0(2) \; ... \;  f_0(n)\\
&f_1(0),f_1(1) \; ... \; f_1(n-1)\\
&...\\
&f_n(0)
\end{split}
\end{align}
with $f_m$ defined recursively as follows:
\begin{equation}\begin{gathered}\label{recursiveF}
f_0(x) = x+1 \\
\forall  m \in \mathbb{N} \quad f_{m+1}(x) = \sum\limits_{a=0}^{x} f_{m}(a)\end{gathered}\end{equation}

Here $A_n$ represents the number in the original stream that would eventually become $M(0,n)$ after \textit{Moessner's Sieve}, and $B_n$ describes the addition on $A_n$ during the  \textit{Moessner's Sieve}: every $m$-th row (starting from $0$-th and ending at $k$-th) denotes the list of $(n-m)$ numbers to be added to $A_n$ at the $m$-th partial sum. (The partial sum before any removal- then-new-stream-formation process is referred to as the $0$-th partial sum.)

Now let us shift our attention to $M(0,n+1)$, which is equivalent to the sum of $A_{n+1}\coloneqq A_n+1$ and $B_{n+1}$ where $B_{n+1}$ is the numbers below:
\begin{align}
\begin{split}
&f_0(0),f_0(1),f_0(2) \; ... \; f_0(n+1)\\
&f_1(0),f_1(1) \; ... \; f_1(n)\\
&...\\
&f_{n+1}(0)
\end{split}
\end{align}
We observe that $B_n$ and $B_{n+1}$ differ in that $B_{n+1}$ has an extra diagonal colored in red as shown below:
\begin{align}
\begin{split}\label{f-extra-d}
&f_0(0),f_0(1),f_0(2) \; ... \;  f_0(n), {\color{red} f_0(n+1)}\\
&f_1(0),f_1(1) \; ... \; f_1(n-1), {\color{red} f_1(n) }\\
&...\\
&f_{n}(0), {\color{red}f_n(1)}\\
&{\color{red} f_{n+1}(0) }
\end{split}
\end{align}
For any $n$, this extra diagonal is equivalent to $V_n$ where
\begin{align}
V_n = \sum\limits_{a=0}^{n+1} f_{a}(n+1-a)
\end{align}
and by (\ref{triangle1}) and (\ref{recursiveF})
\begin{align}\label{V}
V_n = f_0(n+1)+B_{n}
\end{align}
Considering that $f_0(n+1) + 1 = A_{n}$, we have
\begin{align}
\begin{split}
M(0,n+1)-M(0,n) &=  (A_{n+1}-A_{n}) + (B_{n+1}-B_{n}) \\
&= 1 + V_n \\
&= A_n + B_n \\
&= M(0,n)
\end{split}
\end{align}

By demonstrating that $M(0,0)=2^2$ and $\forall n \in \mathbb{N}\; M(0,n+1)= M(0,n) + M(0,n)$ hold, we have inductively showed that
$$\forall n \in \mathbb{N}\; M(0,n)=2^{n+2}$$
as a result of the definitions of multiplication (i.e., $x + x = 2x$) and exponentiation (i.e. $x^n \cdot x = x^{n+1}$).

\section{An inductive proof for \large{$\forall k \in \mathbb{N}\; \small{M(k,0)=(k+2)^n}$}}
The base case (k=0) is the same as that in the previous induction. We are left with demonstrating that $\forall k \in \mathbb{N} \; \phi_M(k,0) \Rightarrow \phi_M(k+1,0)$.

After removing every 2nd element in the original stream, we are left with a stream of odd numbers i.e.,
$$ 1,3,5,7,9 ... $$
Every number in the stream starting from the 3rd (i.e., 5) can be expressed as $2 (k+2) + 1$, with $k$ starts from 0. After the partial sum, we would obtain a new stream
$$1, M(0,0), M(1,0), M(2,0), M(3,0) ...$$
where
\begin{align}\label{n=0}
\forall k \in \mathbb{N}\quad M(k+1,0) = M(k,0) + 2(k + 2) + 1
\end{align}
Therefore we can conclude that
$$\forall k \in \mathbb{N}\; M(k,0)=(k+2)^{2}$$
as a result of the simple fact that $(a+b)^2 = a^2 + 2ab + b^2$.

\section{Our main dish of the day}
If we lay out all instances of $\phi_M(n,k)$ nicely on a piece of paper we would have
\begin{gather*}\color{blue} \phi_M(0,0) \; \phi_M(1,0) \; \phi_M(2,0) \; \phi_M(3,0) \; \phi_M(4,0) \; ... \\
{\color{blue} \phi_M(0,1)} \; \phi_M(1,1) \; \phi_M(2,1) \; \phi_M(3,1) \; \phi_M(4,1) \; ... \\
{\color{blue} \phi_M(0,2)} \; \phi_M(1,2) \; \phi_M(2,2) \; \phi_M(3,2) \; \phi_M(4,2) \; ... \\
{\color{blue} \phi_M(0,3)} \; \phi_M(1,3) \; \phi_M(2,3) \; \phi_M(3,3) \; \phi_M(4,3) \; ... \\
{\color{blue} \phi_M(0,4)} \; \phi_M(1,4) \; \phi_M(2,4) \; \phi_M(3,4) \; \phi_M(4,4) \; ... \\
......
\end{gather*}
where the instances colored in blue have been proven. By demonstrating that
$$\phi_M(k,n) \land \phi_M(k+1,n) \land \phi_M(k,n+1) \Rightarrow \phi_M(k+1,n+1)$$
i.e., for any predicate colored in black above, it is true if its top, left, and top-left neighboring predicates are true, we would cover the entire space and inductively prove that $\phi_M(n,k)$ holds for all $n$ and $k$ in $\mathbb{N}$.

We start off by making the observation that, by \textit{Moessner's Sieve}, $M(k,n)$ is equivalent to the sum of $A_n^k \coloneqq (k+1)(n+2)+1$ and $B_n^k \coloneqq \sum\limits_{i=0}^k \Delta_n^i $ where $\Delta_n^i$ is the sum of the numbers below
\begin{align}
\begin{split}\label{triangle1gen}
&g_0^i(0,n),g_0^i(1,n),g_0^i(2,n) \; ... \;  g_0^i(n,n)\\
&g_1^i(0,n),g_1^i(1,n) \; ... \; g_1^i(n-1,n)\\
&...\\
&g_n^i(0,n)
\end{split}
\end{align}
with $g^i_m$ being a generalization of $f_m$ in (\ref{triangle1}), defined recursively as follows
\begin{equation}\begin{gathered}\label{recursiveFgen}
g_0^i(x,n) = i\cdot(n+2) + x + 1 \\
\forall  m \in \mathbb{N} \quad g_{m+1}^i(x,n) = (\sum\limits_{a=0}^{x} g^i_{m}(a,n))+ \sum\limits_{j=0}^{i-1} \sum\limits_{a=0}^{n-m} g^j_m(a,n) \end{gathered}\end{equation}
Here $A_n^k$ and $B_n^k$ have the same respective interpretations just as $A_n$ and $B_n$ above.  (Note that $g_m^0(x,n) =  f_m(x)$ and that $\sum\limits_{j=0}^{i-1} \sum\limits_{a=0}^{n-m} g^j_m(a,n)$ is simply the sum of the numbers in the first $(n-m)$ rows in all the $\Delta_n^j$ that come before $\Delta_n^i$. In the case when $i$ is 0, there wouldn't be any $\Delta_n^j$, so we don't see the addition of this guy in (\ref{recursiveF}) for the definition of $f_m$.)

For any $n$ and $k$ in $\mathbb{N}$, by  \textit{Moessner's Sieve}, it is clear that
\begin{align}\label{k+1dif}
\begin{split}
M(k+1,n)-M(k,n) &=  (A_{n}^{k+1}-A_{n}^k) + (B_{n}^{k+1}-B_{n}^k) \\
&= (n+2) + \Delta_n^{k+1}
\end{split}
\end{align}
It is also clear that, if $\phi_M(k,n)$ and $\phi_M(k,n+1)$ are true, we have
\begin{align}\label{n+1dif}
\begin{split}
M(k,n+1)-M(k,n) &=  (k+2)^{n+3} - (k+2)^{n+2} \\
&=  M(k,n) \cdot (k + 1) \\
\end{split}
\end{align}
Assuming $\phi_M(k,n)$ and $\phi_M(k,n+1)$, by (\ref{n+1dif}) we have $M(k,n) = M(k,n+1) - M(k,n) \cdot (k + 1)$, and by  (\ref{k+1dif}), we can see that
\begin{align}
\begin{split}\label{stepping-stone2}
&M(k+1,n+1) - M(k,n) \\
=&\;((n+3) + \Delta_{n+1}^{k+1}) + M(k,n) \cdot (k + 1)
\end{split}
\end{align}
On the other hand, assuming $\phi_M(k+1,n)$, by (\ref{n+1dif}) and (\ref{k+1dif}), if
\begin{align}
\begin{split}\label{stepping-stone}
&M(k+1,n+1) - M(k,n) \\
=&\;((n+2) + \Delta_n^{k+1}) + M(k+1,n) \cdot (k+2)
\end{split}
\end{align}
holds, $ M(k+1,n+1) = M(k+1,n) + M(k+1,n) \cdot (k + 2))$ holds, i.e., $\phi_M(k+1,n+1)$ is true.  We now present a proof for
\begin{align}\label{last_inductive_case}
\begin{split}
&((n+3) + \Delta_{n+1}^{k+1}) + M(k,n) \cdot (k + 1)   \\
=\;& ((n+2) + \Delta_n^{k+1}) + M(k+1,n) \cdot (k+2)
\end{split}
\end{align}
assuming only $\phi_M(k,n)$, $\phi_M(k,n+1)$ and $\phi_M(k+1,n)$.

We start off by simplifying the expression into
\begin{align}\label{delta-diff}
\begin{split}
 \Delta_{n+1}^{k+1} - \Delta_{n}^{k+1} = M(k+1,n) \cdot (k+2) - M(k,n) \cdot (k + 1) - 1
\end{split}
\end{align}
and we make the observation that $ \Delta_{n+1}^{k+1}$ and $\Delta_{n}^{k+1}$ differ in that every $g_m^{k+1}$ takes in $n+1$ instead of $n$ in the second parameter, and that there is an extra diagonal (similar to (\ref{f-extra-d})) in $ \Delta_{n+1}^{k+1}$:
\begin{align}
\begin{split}\label{g-extra-d}
&g^{k+1}_0(0,n \;{\color{magenta}+\;1}),g^{k+1}_0(1,n \;{\color{magenta}+\;1}),g^{k+1}_0(2,n \;{\color{magenta}+\;1}) \; ... \;  g^{k+1}_0(n,n \;{\color{magenta}+\;1}), {\color{red} g^{k+1}_0(n+1,n + 1)}\\
&g^{k+1}_1(0,n \;{\color{magenta}+\;1}),g^{k+1}_1(1,n \;{\color{magenta}+\;1}) \; ... \; g^{k+1}_1(n-1,n \;{\color{magenta}+\;1}), {\color{red} g^{k+1}_1(n,n + 1) }\\
&...\\
&g^{k+1}_{n}(0,n \;{\color{magenta}+\;1}), {\color{red}g^{k+1}_n(1,n + 1 )}\\
&{\color{red} g^{k+1}_{n+1}(0,n + 1) }
\end{split}
\end{align}
The difference between $ \Delta_{n+1}^{k+1}$ and $\Delta_{n}^{k+1}$ is thus, by the  \textit{Moessner's Sieve}, the sum of $V_{n+1}^{k+1}$ and $H_{n}^{k+1}$, where $V_{n+1}^{k+1}$ accounts for the extra diagonal (with $n$ being the second parameter in $g_m^{k+1}$), and $H_{n+1}^{k+1}$ accounts for the differences between all $g^{k+1}_m(a,n+1)$ and $g^{k+1}_m(a,n)$ (including those inside the diagonal):
\begin{align}
\begin{split}
V_{n}^{k+1} = \sum\limits_{a=0}^{n+1} g_{a}^{k+1}(n+1-a,n)
\end{split}
\end{align}
and $H_{n}^{k+1}$ is the sum of the numbers below,
\begin{align}
\begin{split}\label{h}
&h^{k+1}_0(0,n),h^{k+1}_0(1,n),h^{k+1}_0(2,n) \; ... \;  h^{k+1}_0(n+1,n) \\
&h^{k+1}_1(0,n),h^{k+1}_1(1,n) \; ... \; h^{k+1}_1(n,n)\\
&...\\
&h^{k+1}_{n+1}(0,n)
\end{split}
\end{align}
where
\begin{align}
\begin{split}
h_m^{i} = g_m^i(x,n+1) - g_m^i(x,n)
\end{split}
\end{align}

Similar to (\ref{V}), by (\ref{triangle1gen}) and (\ref{recursiveFgen}),
\begin{align}
\begin{split}
V_{n}^{k+1} = g_0^{k+1}(n+1,n) + \sum\limits_{i=0}^{k+1} \Delta_n^i
\end{split}
\end{align}
Considering that $g_0^{k+1}(n+1,n) +1 = A_n^{k+1}$, we have
\begin{align}
\begin{split}
V_{n}^{k+1} = M(k+1,n) -1
\end{split}
\end{align}
Since $ \Delta_{n+1}^{k+1} - \Delta_{n}^{k+1} = V_{n}^{k+1} + H_{n}^{k+1}$, (\ref{delta-diff}) can be further reduced into
\begin{align}\label{delta-diff-2}
\begin{split}
H_{n}^{k+1} &=(M^{k+1}_n-M^k_n) \cdot (k+1)  \\
&=((n+2) + \Delta_n^{k+1} ) \cdot (k+1)
\end{split}
\end{align}

It is easy to see that $g_0^i(x,n+1) - g_0^i(x,n) = i$ holds for any $x, n, i \in \mathbb{N}$, so the sum of the first row of numbers in $H_{n}^{k+1}$ is $(n+1)\cdot(k+1)$. We now demonstrate that the bottom $n$ rows sums up to $\Delta_{n}^{k+1} \cdot (k+1)$.
We observe that, for any  $x, n, i, m \in \mathbb{N}$,
\begin{align}
\begin{split}
h_{m+1}^{i}(x,n) = g_m^{i}(x,n) \cdot i
\end{split}
\end{align}
holds, and we give an inductive proof for it. Let $m=0$ be the base case:
\begin{align}
\begin{split}
g_{1}^{i}(x,n+1)-g_{1}^{i}(x,n) = &\text{\{by the defintion of $g$ \}}  \;(\sum\limits_{a=0}^{x} i)+ (\sum\limits_{j=0}^{i-1} (\sum\limits_{a=0}^{n+1} j) + g^j_0(n+1,n)) \\
=&\;(x+1)\cdot i+ (\sum\limits_{j=0}^{i-1} (j\cdot(n+2)) + (j\cdot(n+2) + (n+2))) \\
=&\;(x+1)\cdot i+ (\sum\limits_{j=0}^{i-1} (n+2)\cdot (2\cdot j+1)) \\
= &\;(x+1)\cdot i+  ((n+2)\cdot i \cdot \frac{2 \cdot i+1-1}{2})  \\
= &\;(x+1)\cdot i+ ((n+2) \cdot i \cdot i) \\
=&\; ((x + 1) + (n+2) \cdot i) \cdot i\\
=&\;  g_0^{i}(x,n) \cdot i
\end{split}
\end{align}
The base case is shown to hold. All is left is to demonstrate that $(h_{m+1}^{i}(x,n) = g_m^{i}(x,n) \cdot i) \Rightarrow (h_{m+2}^{i}(x,n) = g_{m+1}^{i}(x,n) \cdot i) $:
\begin{align}
\begin{split}
&g_{m+2}^{i}(x,n+1)-g_{m+2}^{i}(x,n)\\
=&\;(\sum\limits_{a=0}^{x} (g^i_{m+1}(a,n+1)-g^i_{m+1}(a,n)))+ (\sum\limits_{j=0}^{i-1} (\sum\limits_{a=0}^{n-m} g^j_{m+1}(a,n+1) - g^j_{m+1}(a,n)) +  g^j_{m+1}(n-m,n)) \\
&\text{\{by the induction hypothesis\}} \\
=&\;(\sum\limits_{a=0}^{x}  g_m^{i}(a,n) \cdot i)+ (\sum\limits_{j=0}^{i-1} (\sum\limits_{a=0}^{n-m} g_m^{j}(a,n) \cdot j) +  g^j_{m+1}(n-m,n)) \\
=&\;(\sum\limits_{a=0}^{x}  g_m^{i}(a,n) \cdot i)+ (\sum\limits_{j=0}^{i-1} (\sum\limits_{a=0}^{n-m} g_m^{j}(a,n) \cdot j)) +  (\sum\limits_{j=0}^{i-1} g^j_{m+1}(n-m,n) ) \\
=&\;(\sum\limits_{a=0}^{x}  g_m^{i}(a,n) \cdot i)+ (\sum\limits_{j=0}^{i-1} (\sum\limits_{a=0}^{n-m} g_m^{j}(a,n) \cdot j)) +  (\sum\limits_{j=0}^{i-1}(\sum\limits_{a=0}^{n-m} g_m^{j}(a,n) \cdot (i-j) )) \\
=&\;(\sum\limits_{a=0}^{x} g^i_{m}(a,n) \cdot i) + (\sum\limits_{j=0}^{i-1} \sum\limits_{a=0}^{n-m} g^j_m(a,n)\cdot i) \\
=&\; g_{m+1}^{i}(x,n) \cdot i
\end{split}
\end{align}
where $\sum\limits_{j=0}^{i-1} g^j_{m+1}(n-m,n) = \sum\limits_{j=0}^{i-1} \sum\limits_{a=0}^{n-m} g_m^{j}(a,n) \cdot (i-j)$ can be shown as follows
\begin{align}
\begin{split}
\sum\limits_{j=0}^{i-1} g^j_{m+1}(n-m,n) =& \text{\{by the defintion of $g$ \}} \sum\limits_{j=0}^{i-1} (\sum\limits_{a=0}^{n-m} g^j_{m}(a,n)+ \sum\limits_{u=0}^{j-1} \sum\limits_{a=0}^{n-m} g^u_m(a,n)) \\
=& \sum\limits_{j=0}^{i-1} \sum\limits_{a=0}^{n-m} (g^j_{m}(a,n)+ \sum\limits_{u=0}^{j-1} g^u_m(a,n)) \\
=& \sum\limits_{a=0}^{n-m} \sum\limits_{j=0}^{i-1}  \sum\limits_{u=0}^{j} g^u_m(a,n)) \\
=& \sum\limits_{a=0}^{n-m} \sum\limits_{u=0}^{i-1}  \sum\limits_{j=u}^{i-1} g^u_m(a,n)) \\
=&\sum\limits_{a=0}^{n-m}  \sum\limits_{u=0}^{i-1} g_m^{u}(a,n) \cdot (i-u)
\end{split}
\end{align}
Therefore, $h_{m+1}^{i}(x,n) = g_m^{i}(x,n) \cdot i$ holds and we obtain
\begin{align}
\begin{split}
H_{n}^{k+1} =(n+2) \cdot (k+1) + \Delta_n^{k+1}  \cdot (k+1)
\end{split}
\end{align}
proving that (\ref{delta-diff-2}) holds, and therefore (\ref{last_inductive_case}) holds, consequently indicating that
$$M(k+1,n+1)-M(k+1,n)  = M(k+1,n) \cdot (k + 2) $$
holds assuming only $\phi_M(k,n)$, $\phi_M(k,n+1)$ and $\phi_M(k+1,n)$, i.e.,
$$\forall n,k \in \mathbb{N}\quad \phi_M(k,n) \land \phi_M(k+1,n) \land \phi_M(k,n+1) \Rightarrow \phi_M(k+1,n+1)$$
Q.E.D


\begin{thebibliography}{9}
\bibitem{moessner}
Moessner, A.
\emph{Eine Bemerkung �ber die Potenzen der nat�rlichen Zahlen}.
Aus den Sitzungsberichten
der Bayerische Akademie der Wissenschaften, Mathematisch-naturwissenschaftliche Klasse 1951 Nr. 3
(1951)
\bibitem{perron}
Perron,  O.
\emph{Beweis des Moessnerschen Satz}.
Aus den Sitzungsberichten der Bayerische Akademie der
Wissenschaften, Mathematisch-naturwissenschaftliche Klasse 1951 Nr. 4, (1951)
\bibitem{p}
Paasche, I.
\emph{Ein neuer Beweis des moessnerischen Satzes.}
Sitzungsberichten der Bayerischen Akademie der
Wissenschaften, Mathematischnaturwissenschaftliche Klasse 1952, 1:1-5 (1953).
\bibitem{l}
Long, C. T.
\emph{On the Moessner's theorem on integral powers}.
Amer. Math. Monthly, 73(8):846-851, October (1966)
\bibitem{h}
Hinze, R.
\emph{Scans and convolutions - a calculational proof of Moessner's theorem.}
In Post-proceedings of the 20th International Symposium on the Implementation and Application of Functional
Languages (IFL '08), volume 5836 of Lecture Notes in Computer Science (2009)
\bibitem{n_r}
Niqui, M. and Rutten, J.J.M.M.
\emph{An exercise in coinduction: Moessner's theorem.}
Technical Report SEN-1103, CWI Amsterdam (2011)
\bibitem{k_s}
Kozen, D. and Silva, A.
\emph{On Moessner's Theorem}. (2011)
\end{thebibliography}
\end{document}